\documentclass[11pt]{amsart}
\usepackage{mathrsfs}
\usepackage{amssymb}

\usepackage{titletoc}
\usepackage{amscd}
\usepackage{amsmath, amssymb}
\usepackage{amsfonts}
\usepackage[colorlinks,linkcolor=blue,citecolor=blue, pdfstartview=FitH]{hyperref}

\numberwithin{equation}{section}

\theoremstyle{plain}
\newtheorem{thm}{Theorem}[section]
\newtheorem{theorem}[thm]{Theorem}
\newtheorem{lemma}[thm]{Lemma}
\newtheorem{corollary}[thm]{Corollary}
\newtheorem{proposition}[thm]{Proposition}

\theoremstyle{definition}

\newtheorem{problem}[thm]{Problem}
\newtheorem{remark}[thm]{Remark}

\newtheorem{definition}[thm]{Definition}

\newtheorem{example}[thm]{Example}

\newtheorem{defn-thm}[thm]{Definition-Theorem}

\newcommand{\sO}{{\mathcal O}}

\newcommand{\C}{{\mathbb C}}

\newcommand{\Q}{{\mathbb Q}}
\newcommand{\R}{{\mathbb R}}

\newcommand{\qtq}[1]{\quad\mbox{#1}\quad}
\newcommand{\bp}{\bar{\partial}}
\newcommand{\Om}{\Omega}

\newcommand{\ts}{\otimes}

\newcommand{\btheorem}{\begin{theorem}}
\newcommand{\etheorem}{\end{theorem}}
\newcommand{\bproposition}{\begin{proposition}}
\newcommand{\eproposition}{\end{proposition}}
\newcommand{\bdefinition}{\begin{definition}}
\newcommand{\edefinition}{\end{definition}}
\newcommand{\bcorollary}{\begin{corollary}}
\newcommand{\ecorollary}{\end{corollary}}
\newcommand{\bproof}{\begin{proof}}
\newcommand{\eproof}{\end{proof}}
\newcommand{\bremark}{\begin{remark}}
\newcommand{\eremark}{\end{remark}}
\newcommand{\eexample}{\end{example}}
\newcommand{\bexample}{\begin{example}}

\newcommand{\elemma}{\end{lemma}}
\newcommand{\blemma}{\begin{lemma}}

\newcommand{\sq}{\sqrt{-1}}

\newcommand{\p}{\partial}

\renewcommand{\bar}{\overline}
\newcommand{\eps}{\varepsilon}

\renewcommand{\phi}{\varphi}

\newcommand{\ee}{\end{eqnarray*}}
\newcommand{\be}{\begin{eqnarray*}}

\newcommand{\beq}{\begin{equation}}
\newcommand{\eeq}{\end{equation}}

\newcommand{\bd}{\begin{enumerate}}
\newcommand{\ed}{\end{enumerate}}

\renewcommand{\tilde}{\widetilde}

\renewcommand{\>}{\rightarrow}

\usepackage{fancyhdr}
\begin{document}
\title{A partial
converse to the Andreotti-Grauert theorem } \makeatletter
\let\uppercasenonmath\@gobble
\let\MakeUppercase\relax
\let\scshape\relax
\makeatother
\author{Xiaokui Yang}
\date{}
\address{Morningside Center of Mathematics, Academy of Mathematics and\\ Systems Science, Chinese Academy of Sciences, Beijing, 100190, China}
\address{HCMS, CEMS, NCNIS, HLM, UCAS, Academy of Mathematics and\\ Systems Science, Chinese Academy of Sciences, Beijing 100190,
China} \email{\href{mailto:xkyang@amss.ac.cn}{{xkyang@amss.ac.cn}}}

\maketitle

\maketitle

\begin{abstract}  Let $X$ be a smooth
projective manifold  with $\dim_\C X=n$. We show that if
 a line bundle $L$ is $(n-1)$-ample, then it is $(n-1)$-positive. This is a
partial converse to the Andreotti-Grauert theorem. As an
application, we show that a projective manifold $X$ is uniruled if
and only if there exists a Hermitian metric $\omega$ on $X$ such
that its Ricci curvature $\mathrm{Ric}(\omega)$ has at least one
positive eigenvalue everywhere.

\end{abstract}
\setcounter{tocdepth}{1} \tableofcontents

\section{Introduction}

One of the most fundamental topics in algebraic geometry is  the
characterizing  ampleness of  line bundles by using numerical and
cohomological vanishing theorems. The theorem of
Cartan--Serre--Grothendieck is a milestone in this direction: the
following statements are equivalent \bd
 \item  $L$ is an ample line bundle over a projective manifold $X$;
\item For every coherent sheaf $\mathcal F$ on $X$, there exists a positive integer $m_0=m_0(X,\mathcal F,L)$ such that
$$
 H^{i}(X,\mathcal F\otimes L^{\otimes m})=0
$$
for all $i>0$ and all $m\geq m_0$. \ed On complex projective
manifolds, ampleness is also equivalent to the existence of a smooth
metric with positive curvature, thanks to the celebrated Kodaira
embedding theorem.

In \cite{AG62}, Andreotti and Grauert considered the case of line
bundles with curvature of mixed signature. Now it is formulated into

\begin{definition}\label{Def} Let $L$ be a holomorphic line bundle over a  compact complex
manifold $X$ with $\dim_\C X=n$. \bd \item  $L$ is called
\emph{$q$-positive}, if there exists a smooth Hermitian metric $h$
on $L$ such that the Chern curvature $R^{(L,h)}=-\sq\p\bp\log h$ has
at least $(n-q)$ positive eigenvalues at every point on $X$. \item
$L$ is called \textit{ $q$-ample}, if for any coherent sheaf
$\mathcal{F}$ on $X$ there exists a positive integer
$m_{0}=m_{0}(X,L,\mathcal{F})>0$ such that \beq H^{i}\left( X,
\mathcal{F}\otimes L^{ m}\right) =0\ \ \ \mathrm{for}\  i>q,\ m\geq
m_{0}. \label{fv}\eeq \ed
\end{definition}

\noindent It is obvious that $L$ is $0$-positive if and only if $L$
is positive, and $L$ is $0$-ample if and only if $L$ is ample. Hence
the $0$-positivity and $0$-ampleness are equivalent. In
\cite[Theorem~14]{AG62}, Andreotti and Grauert proved the following
fundamental theorem (see also \cite[Proposition~2.1]{DPS96}).

\btheorem[Andreotti-Grauert]\label{C} A $q$-positive line bundle is
$q$-ample. \etheorem

\noindent Historically, the Andreotti-Grauert Theorem is the first
result on the relationship between (partially) positive line bundles
and the cohomological vanishing theorems. In the pioneer work
\cite{DPS96}, Demailly, Peternell and Schneider systematically
investigated partial vanishing theorems and proposed the following
problem on the converse to the Andreotti-Grauert theorem:

\begin{problem}\label{Q} On a projective manifold,  if a line bundle is $q$-ample, is it $q$-positive?
\end{problem}

\noindent This is a long-standing open problem. The key difficulty
arises from constructing a precise metric
 according to the formal  partial vanishing theorem (\ref{fv}). Recently, there are some progress on this problem, mainly contributed by   Demailly, Totaro, Ottem, K$\mathrm{\ddot{u}}$ronya,
 Matsumura, Brown and etc (see \cite{Dem11, Tot10, Mat12, Ott12,
Bro12, Kur10, GK15} and also the references therein).  Totaro proved
in \cite{Tot10} that the notion of $q$-ampleness is equivalent to
others previously studied in \cite{DPS96}. As a result, the
$q$-ampleness of a line bundle depends only on its numerical class,
and the cone of such bundles is open. In particular, Totaro
established that the $(n-1)$-ample cone of an $n$-dimensional
projective manifold $X$ is equal to the negative of the complement
of the pseudo-effective cone of $X$ (see also Corollary \ref{D}).
 In dimension
two, Demailly
 proved in \cite{Dem11} an asymptotic version of this converse to the
Andreotti-Grauert Theorem using tools related to the holomorphic
Morse inequality and asymptotic cohomology; subsequently,  Matsumura
obtained in \cite[Theorem~1.3]{Mat12} a positive answer to the
question for projective surfaces. However, there exist
higher-dimensional counterexamples to the converse
Andreotti--Grauert problem in the range $\tfrac{\dim X}{2}-1<q<\dim
X-2$, constructed by Ottem \cite[Theorem~10.3]{Ott12}.
 Our main result in this paper is a partial converse to the
Andreotti-Grauert theorem on smooth projective manifolds.

\btheorem\label{A}  Let $L$ be a line bundle over a  smooth
projective manifold $X$ with $\dim_\C X=n$. If $L$ is $(n-1)$-ample,
then it is $(n-1)$-positive. \etheorem

\noindent In particular, when $\dim_\C X=2$, the converse
Andreotti-Grauert Problem \ref{Q} is true (see also
\cite[Theorem~1.3]{Mat12}). Actually, Theorem \ref{A} is a
straightforward application of the following result on general
\textit{compact complex manifolds}.

\btheorem\label{B} Let $X$ be a  compact complex manifold with
$\dim_\C X=n$. Then the following statements are equivalent:

\bd

\item  $L$ is $(n-1)$-positive;

\item  the dual line bundle $L^{-1}$ is not pseudo-effective.

\ed \etheorem \noindent Note that, Theorem \ref{B} is also valid if
we replace the line bundle $L$ by a Bott-Chern  class $\alpha \in
H^{1,1}_{\mathrm{BC}}(X)$ (see Theorem \ref{B2}). The key
ingredients in the proof of Theorem \ref{B} rely on several results
in our previous paper \cite{Yang17} on geometric characterizations
of pseudo-effective line bundles (resp. Bott-Chern classes).
 As an application of Theorem \ref{A}, Theorem \ref{B} and
Theorem \ref{C}, one has

\bcorollary\label{D} On a projective manifold $X$ of complex
dimension $n$, the following are equivalent:

\bd \item $L$ is $(n-1)$-ample;

\item $L$ is $(n-1)$-positive;

\item $L^{-1}$ is not pseudo-effective.

\ed

\ecorollary \noindent   Note  that the equivalence of $(1)$ and
$(3)$ is also obtained in \cite[Proposition~2.5]{DPS96},
\cite[Theorem~9.1]{Tot10} and \cite[Lemma~2.4]{Mat12} by different
methods in one direction.

According to Ottem's counter-examples in \cite[Theorem~10.3]{Ott12}, for $\tfrac{\dim X}{2}-1<q<\dim
X-2$, the $q$-ampleness can not imply the $q$-positivity. On the contrary,
by Theorem \ref{A}, we obtain

\bproposition \label{F} Let $X$ be a smooth projective manifold with $\dim_\C X=n$. Suppose   $L$ is $q$-ample, then
 the restriction of $L$ to every  codimension-$(n-q-1)$ smooth submanifold $Y$ is $q$-positive.
\eproposition

\noindent In particular, we have
\bcorollary Let $X$ be a smooth projective manifold with $\dim_\C X=n$. If $L$ is $(n-2)$-ample, then
 the restriction of $L$ to every  codimension-$1$ smooth submanifold is $(n-2)$-positive.

\ecorollary

\noindent On the other hand, by using the classical result of
\cite{BDPP13} and Yau's theorem \cite{Yau78}, one obtains

\bcorollary\label{E} On a projective manifold $X$ of complex
dimension $n$, the following are equivalent:

\bd \item $X$ is uniruled;

\item $K_X^{-1}$ is $(n-1)$-positive;

\item  there exists a smooth Hermitian metric $\omega$ on $X$ such that
its Ricci curvature $\mathrm{Ric}(\omega)$ has at least one positive
eigenvalue everywhere;

\item $K_X$ is not pseudo-effective;

\item $K^{-1}_X$ is $(n-1)$-ample.

 \ed

\ecorollary

\vskip 1\baselineskip

\bremark The equivalence of $(1)$ and $(3)$ is also conjectured by
S.-T Yau.

\eremark

\bremark From Theorem \ref{A}, Theorem \ref{B}, Corollary \ref{D},
Corollary \ref{E}, Theorem \ref{B1} and Theorem \ref{B2}, one can
derive and formulate various cone dualities on compact complex
manifolds. The vector bundle analogous of the main results and
applications are  obtained in \cite{Yang18}. \eremark

\noindent\textbf{Acknowledgement.}  I am very grateful to Professor
 Kefeng Liu for his support, encouragement and stimulating
discussions over  years. I would also like to thank Professors
Junyan Cao, J.-P. Demailly, S. Matsumura, Yum-Tong Siu, Junchao
Shentu, Xiaotao Sun, Valentino Tosatti, Jian Xiao and Xiangyu Zhou
for  some useful suggestions.  The author would also like to thank
the anonymous referees whose comments and suggestions helped improve
and clarify the paper. This work was partially supported   by
China's Recruitment
 Program of Global Experts and  NSFC 11688101.

\vskip 3\baselineskip

\section{Partially positive line bundles on compact complex manifolds}

In this section, we deal with two different notions on $q$-positive
line bundles over compact complex manifolds.

\begin{definition} Let X be a compact complex manifold and $L$ be a holomorphic line
bundle over $X$. \bd \item  $L$ is called \emph{$q$-positive}, if
there exists a smooth Hermitian metric $h$ on $L$ and a smooth
Hermitian metric $\omega$ on $X$ such that the Chern curvature
$R^{(L,h)}=-\sq\p\bp\log h$ has at least $(n-q)$ positive
eigenvalues at any point on $X$ (with respect to $\omega$);
\item $L$ is called \textit{ uniformly $q$-positive}, if there
exists a Hermitian metric $h$ on $L$ nd a smooth Hermitian metric
$\omega$ on $X$
 such that  the summation
of any distinct $(q+1)$ eigenvalues (counting multiplicity) of the
Chern curvature $R^{(L,h)}=-\sq\p\bp\log h$ is positive at any point
of $X$ (with respect to $\omega$). \ed

\end{definition}

\noindent The following result  is obtained by changing the metric
on the complex manifold and the background idea dates back to
\cite[Section~5]{AV65}.

\bproposition\label{equivalent}  The following statements are
equivalent: \bd
\item $L$ is $q$-positive;
\item $L$ is uniformly $q$-positive.
\ed \eproposition

\bproof $(2)\Longrightarrow(1)$. Let $\omega$ be a Hermitian metric
on $X$ and $h$ be a smooth Hermitian metric on $L$ such that
$R^{(L,h)}=-\sq\p\bp\log h$ is uniformly $q$-positive. Let
$\lambda_1\geq \cdots\geq \lambda_n$ be the eigenvalues of
$R^{(L,h)}$ with respect to $\omega$ over some coordinate chart. We
have $\lambda_{n-q}>0$. Otherwise, the summation of $q+1$
eigenvalues $\lambda_{n-q}+\lambda_{n-q+1}+\cdots+\lambda_n \leq 0$.

$(1)\Longrightarrow (2).$ We assume that there exists a smooth
Hermitian metric $h$ on
 $L$ such that the curvature $R^{(L,h)}=-\sq\p\bp\log h$ has at least
 $(n-q)$ positive eigenvalues at each point $p\in X$. Let $\omega_0$ be a
 fixed Hermitian metric on $X$. For simplicity, we denote by $R$ and
 $\Om$ the local matrix representations of the matrices $R^{(L,h)}$
 and $\omega_0$ respectively, in some local holomorphic frames of $X$. Let
 $$\lambda_1(z)\geq \cdots\geq \lambda_n(z)$$
be the eigenvalues of $R^{(L,h)}$ with respect to $\omega_0$. It is
obvious that $\lambda_1,\cdots, \lambda_n$ are eigenvalues of the
matrix $R\Om^{-1}$. Note that $R\Om^{-1}$  represents a tensor in
$\Gamma(X,\mathrm{End}(T^{1,0}X))$, and so the eigenvalues
$\lambda_i$ are independent of the choice of coordinates.
 Since $\lambda_{n-q}$ is a continuous function, we set
\beq \lambda_0=\frac{\log (n+1)}{\inf_X
\lambda_{n-q}}.\label{para}\eeq $\lambda_0$ is a positive number
since $L$ is $q$-positive and $X$ is compact. We define a new
Hermitian metric $\omega$ over $X$ with local matrix representation
$\tilde \Om$ by the following formula \beq \tilde
\Om^{-1}=\Om^{-1}\cdot \left(Id+
\sum_{k=1}^\infty\frac{\lambda_0^k(R\Om^{-1})^k}{(k+1)!}\right).\eeq
Note that the matrix $\tilde \Om^{-1}$ is positive definite. Indeed,
the eigenvalues of the matrix $Id+
\sum_{k=1}^\infty\frac{\lambda_0^k(R\Om^{-1})^k}{(k+1)!}$ are given
by
$$1+\sum_{k=1}^\infty\frac{\lambda_0^k\lambda_i^k}{(k+1)!}=\frac{e^{\lambda_0\lambda_i}-1}{\lambda_0\lambda_i}> 0,\ \  \text{if}\ \ \  \lambda_i\neq 0.$$
It is not hard to see that the Hermitian metric $\omega$ is globally
well-defined on $X$. Let $\kappa_1\geq \cdots\geq \kappa_n$ be the
eigenvalues of $R^{(L,h)}$ with respect to the new  metric $\omega$,
i.e. they are the eigenvalues of $R\tilde \Om^{-1}$. Note also that
$$R\tilde
\Om^{-1}=\lambda_0^{-1}\left(\sum_{k=0}^\infty\frac{\lambda_0^k(R\Om^{-1})^k}{k!}-Id\right).$$
A straightforward computation shows
$$\kappa_{n-q}=\frac{e^{\lambda_0\lambda_{n-q}}-1}{\lambda_0}\qtq{and}
\kappa_n=\frac{e^{\lambda_0\lambda_n}-1}{\lambda_0}.$$  For any
summation of $(q+1)$ (distinct) eigenvalues of $R^{(L,h)}$ with
respect to the new metric $\omega$, we have the inequality
\be \sum_{\ell=1}^{q+1}\kappa_{i_\ell}&\geq & \kappa_{n-q}+\cdots+\kappa_n\\&\geq& \kappa_{n-q}+q\kappa_n\\
&\geq &
\lambda_0^{-1}\left(e^{\lambda_0\lambda_{n-q}}+qe^{\lambda_0\lambda_n}-(q+1)\right)\\
&> &\lambda_0^{-1}\left(e^{\lambda_0\lambda_{n-q}}-(q+1)\right)>0
\ee since
$e^{\lambda_0\lambda_{n-q}}=e^{\log(n+1)\frac{\lambda_{n-q}}{\inf_X
\lambda_{n-q}}}\geq n+1$ by (\ref{para}). \eproof

\vskip 1\baselineskip

\noindent The following special case of Proposition \ref{equivalent}
is of particular interest in complex geometry.
\bcorollary\label{scalar} The following statements are equivalent:
\bd
\item $L$ is $(n-1)$-positive;
\item there exists a smooth Hermitian metric $h$ on $L$ and a Hermitian metric $\omega$ on
$X$ such that the  function \beq \mathrm{tr}_\omega(-\sq\p\bp\log
h)>0.\eeq
  \ed
\ecorollary

\bremark The function $\mathrm{tr}_\omega(-\sq\p\bp\log h)$ is
globally defined on $X$ and it is independent of the choice of
coordinates. It is also called the scalar curvature of the Chern
curvature $R^{(L,h)}=-\sq\p\bp\log h$ with respect to the Hermitian
metric $\omega$. \eremark

\vskip 2\baselineskip

\section{The pseudo-effective line bundles on compact complex manifolds}

A line bundle $L$ on a compact complex manifold $X$ is called
\emph{pseudo-effective} if there exists a (possibly) singular
Hermitian metric $h$ on $L$ such that its Chern curvature
$R^{(L,h)}=-\sq\p\bp\log h\geq 0$ in the sense of current. In order
to describe pseudo-effective line bundles in a differential
geometric setting, we introduce the Bott-Chern
  cohomology  on $X$：
  $$ H^{p,q}_{\mathrm{BC}}(X):= \frac{\text{Ker} d \cap \Om^{p,q}(X)}{\text{Im} \p\bp \cap \Om^{p,q}(X)}.$$

\noindent Let $\mathrm{Pic}(X)$ be the set of holomorphic line
bundles over $X$. As similar as the first Chern class map
$c_1:\mathrm{Pic}(X)\>H^{1,1}_{\bp}(X)$, there is a \emph{first
Bott-Chern class} map \beq c_1^{\mathrm{BC}}:\mathrm{Pic}(X)\>
H^{1,1}_{\mathrm{BC}}(X).\eeq Given any holomorphic line bundle
$L\to X$ and any Hermitian metric $h$ on $L$, we define
$c_1^{\mathrm{BC}}(L)$ to be the class of its curvature form
$R^{(L,h)}=-\sq\p\bp\log h$ in $H^{1,1}_{\mathrm{BC}}(X)$ (modulo a
constant $2\pi$).  A Hermitian metric $\omega$ is called a
\emph{Gauduchon metric} if $\p\bp\omega^{n-1}=0$ where $\dim_\C
X=n$. It is proved by Gauduchon (\cite{Ga1}) that, in the conformal
class of each Hermitian metric, there exists a unique Gauduchon
metric (up to constant scaling).

\bproposition\label{pseudo} The following statements are equivalent

\bd \item $L$ is  pseudo-effective;

\item For any  Gauduchon metric $\omega_{\mathrm{G}}$ on $X$, one has \beq \int_X c_1^{\mathrm{BC}}(L)\wedge
\omega_{\mathrm{G}}^{n-1}\geq 0.\eeq \ed

\eproposition

\bproof The proof is essentially contained in
\cite[Theorem~1.1]{Yang17} or \cite[Theorem~3.4]{Yang17} which
relies on Lamari's positivity criterion (\cite{La99}) and an
observation of Michelsohn (\cite{Mi}). For readers' convenience, we
include a proof here.

 $(1)\Longrightarrow (2)$.  Suppose $L$ is
pseudo-effective, it is well-known that there exist a smooth
Hermitian metric $h$ on $L$ and a real valued
  function $\phi\in  \mathscr{L}^1(X,\R)$ such that
  $$R^{(L,h)}+\sq\p\bp\phi\geq 0$$
in the sense of current where $R^{(L,h)}=-\sq\p\bp\log h$. Then for
any smooth Gauduchon metric $\omega_{\mathrm{G}}$, we have \be
\mathscr \int_X c_1^{\mathrm{BC}}(L)\wedge
\omega_{\mathrm{G}}^{n-1}&=&
\int_XR^{(L,h)}\wedge \omega_{\mathrm{G}}^{n-1}\\
&=&\left(R^{(L,h)}+\sq\p\bp\phi,\omega_{\mathrm{G}}^{n-1}\right)\geq
0\ee since $\p\bp\omega_{G}^{n-1}=0$ and $R^{(L,h)}+\sq\p\bp\phi\geq
0$ in the sense of current.\\

$(2)\Longrightarrow (1)$. We define several sets:

\bd
\item[$\bullet$] $\mathscr E$ is the set of real $\p\bp$-closed $(n-1,n-1)$
forms on $X$;

\item[$\bullet$] $\mathscr V$ is the set of real positive $\p\bp$-closed $(n-1,n-1)$
forms on $X$;

\item[$\bullet$] $\mathscr G=\{\omega^{n-1}\ |\ \omega\ \text{is a Gauduchon metric
}\}$.

\ed

\noindent In \cite[p.279-p.280]{Mi}, M.L. Michelsohn observed that
$\mathscr V=\mathscr G$. Let $\mathscr W$ be the space of smooth
Gauduchon metrics on  $X$.  We also define $\mathscr F:\mathscr
W\>\R$ by
$$\mathscr F(\omega)=\int_X c_1^{\mathrm{BC}}(L)\wedge \omega^{n-1}.$$
Hence, by the assumption of $(2)$, we have $\mathscr F(\omega)\geq
0$ for every $\omega\in \mathscr W$. Fix an arbitrary smooth
Hermitian metric $h$ on $L$. Since $\mathscr V=\mathscr G$,  for any
$\p\bp$-closed positive $(n-1,n-1)$ form $\psi\in \mathscr V$, there
exists a smooth Gauduchon metric $\omega$ such that
$\omega^{n-1}=\psi$. Hence \beq \int_X R^{(L,h)}\wedge \psi=\int_X
R^{(L,h)}\wedge \omega^{n-1}=\int_X c_1^{\mathrm{BC}}(L)\wedge
\omega^{n-1}=\mathscr F(\omega)\geq 0. \eeq That means, as a
functional on $\mathscr V$,   $R^{(L,h)}$ is non-negative. Note that
$\mathscr V$ is a hyperplane in  $\mathscr E$.  As proved in
\cite[Lemma~3.3]{La99}, by Hahn-Banach theorem,
$c_1^{\mathrm{BC}}(L)$ is pseudo-effective, and there exists a
locally integrable function $\phi\in \mathscr L^1(X,\R)$ such that
$$R^{(L,h)}+\sq\p\bp\phi\geq 0$$ in the sense of
current. That means, $L$ is pseudo-effective. \eproof

\noindent Of course, Proposition \ref{pseudo} has the following
variant:

\bproposition\label{pseudo2} The following statements are equivalent

\bd \item the dual line bundle $L^{-1}$ is not pseudo-effective;

\item there exists a Gauduchon metric $\omega_{\mathrm{G}}$ such
that $$\int_X c_1^{\mathrm{BC}}(L)\wedge
\omega_{\mathrm{G}}^{n-1}>0.$$ \ed\eproposition

\vskip 2\baselineskip

\section{The proofs of Theorem \ref{A} and Theorem \ref{B}}
In this section, we prove Theorem \ref{A}, Theorem \ref{B} and Proposition \ref{F}.\\

\noindent\emph{The proof of Theorem \ref{B}.} $(1)\Longrightarrow
(2).$  By Corollary \ref{scalar}, there exist a smooth Hermitian
metric $h$ on $L$ and a Hermitian metric $\omega$ on $X$ such that
the function \beq \mathrm{tr}_\omega(-\sq\p\bp\log
h)>0.\label{core}\eeq Let $\omega_{\mathrm{G}}=e^f\omega$ be a
Gauduchon metric in the conformal class of $\omega$ (\cite{Ga2}),
then by (\ref{core}), we obtain
$$\mathrm{tr}_{\omega_{\mathrm{G}}} R^{(L,h)}=e^{-f}\cdot \mathrm{tr}_{\omega} R^{(L,h)}>0,$$
where $R^{(L,h)}=-\sq\p\bp\log h$. In particular, we have \beq
\int_X \text{tr}_{\omega_G} R^{(L,h)}\cdot
\omega_{\mathrm{G}}^{n}=n\int_X R^{(L,h)}\wedge \omega_G^{n-1}=
n\int_X c_1^{\mathrm{BC}}(L)\wedge \omega_G^{n-1}>0.\eeq  By
Proposition \ref{pseudo2}, the dual line bundle $L^{-1}$ is not
pseudo-effective.\\

 $(2)\Longrightarrow (1)$. If $L^{-1}$ is not pseudo-effective, by Proposition \ref{pseudo2}, there exists a Gauduchon metric $\omega_{\mathrm{G}}$ such
that $$\int_X c_1^{\mathrm{BC}}(L)\wedge
\omega_{\mathrm{G}}^{n-1}>0.$$ We shall use Gauduchon's conformal
trick (\cite{Ga1, Ga3}, see also \cite{Yang16,
 Yang17}) to construct a smooth Hermitian metric $h$ on $L$ such that $\mathrm{tr}_{\omega_{\mathrm G}}(-\sq\p\bp\log h)>0$. Hence, by Corollary \ref{scalar}, $L$ is $(n-1)$-positive.

  Fix a smooth Hermitian metric $h_0$ on $L$. Let $$R_0=-\sq\p\bp\log
 h_0$$ be the Chern curvature of $(L,h_0)$. It is easy to see that
\beq \int_X \text{tr}_{\omega_G} R_0\cdot
\omega_{\mathrm{G}}^{n}=n\int_X R_0\wedge \omega_G^{n-1}= n\int_X
c_1^{\mathrm{BC}}(L)\wedge \omega_G^{n-1}\eeq Since
$\omega_{\mathrm{G}}$ is Gauduchon, i.e.
$\p\bp\omega^{n-1}_{\mathrm{G}}=0$  and the integration
$$\int_X\left(\text{tr}_{\omega_G} R_0-\frac{n\int_X
c_1^{\mathrm{BC}}(L)\wedge
\omega_G^{n-1}}{\int_X\omega^n_G}\right)\omega_G^n=0,$$ the equation
 \beq \text{tr}_{\omega_G}\sq\p\bp f=\text{tr}_{\omega_G} R_0- \frac{n\int_X
c_1^{\mathrm{BC}}(L)\wedge \omega_G^{n-1}}{\int_X\omega^n_G}\eeq has
a solution $f\in C^\infty(X)$ which is well-known(e.g. \cite{Ga1}，
or \cite[Theorem~2.2]{CTW16}).
 Let $h=e^f\cdot h_0$ be a smooth
Hermitian metric on $L$. The Hermitian line bundle $(L,h)$ has Chern
curvature
$$R^{(L,h)}=-\sq\p\bp\log h=R_0-\sq \p\bp f.$$
The scalar curvature of $R^{(L,h)}$ with respect to
$\omega_{\mathrm{G}}$ is
$$\mathrm{tr}_{\omega_{\mathrm{G}}} R^{(L,h)}=\text{tr}_{\omega_G} R_0-\text{tr}_{\omega_G}\sq\p\bp f=\frac{n\int_X
c_1^{\mathrm{BC}}(L)\wedge \omega_G^{n-1}}{\int_X\omega^n_G}>0.$$
The proof of Theorem \ref{B} is completed. \qed

\vskip 1\baselineskip

\noindent By Corollary \ref{scalar}, Proposition \ref{pseudo2} and
Theorem \ref{B}, we obtain

\btheorem\label{B1} Let $L$ be a line bundle over a compact complex
manifold $X$ with $\dim_\C X=n$. The following statements are
equivalent:

\bd \item The dual line bundle $L^{-1}$ is not pseudo-effective;

 \item There exists a smooth Gauduchon metric $\omega_{\mathrm{G}}$ on
$X$ such that
$$\int_X c_1^{\mathrm{BC}}(L)\cdot \omega_{\mathrm{G}}^{n-1}>0;$$

\item There exist a smooth Hermitian metric $h$ on $L$ and a
Hermitan metric $\omega$ on $X$ such that the scalar curvature of
the Chern curvature $R^{(L,h)}=-\sq\p\bp\log h$ with respect to
$\omega$ is positive, i.e.
$$s=\mathrm{tr}_\omega R^{(L,h)}>0;$$
\item $L$ is $(n-1)$-positive.
\ed \etheorem

\noindent It is easy to see that Theorem \ref{B1} has the following
version on Bott-Chern classes:\btheorem\label{B2} Let $X$ be a
compact complex manifold $X$ with $\dim_\C X=n$ and $[\alpha] \in
H^{1,1}_{\mathrm{BC}}(X)$. The following statements are equivalent:

\bd \item The class $-[\alpha]_{\mathrm{BC}}$ is not
pseudo-effective;

 \item There exists a smooth Gauduchon metric $\omega_{\mathrm{G}}$ on
$X$ such that
$$\int_X [\alpha]_{\mathrm{BC}}\cdot \omega_{\mathrm{G}}^{n-1}>0;$$

\item There exist a smooth $(1,1)$ form  $\chi\in [\alpha]_{\mathrm{BC}}$  and a
Hermitian metric $\omega$ on $X$ such that
$$\mathrm{tr}_\omega \chi>0;$$
\item $[\alpha]_{\mathrm{BC}}$ is $(n-1)$-positive.
\ed \etheorem

\vskip 1\baselineskip

\noindent Now we are ready to prove our main theorem.\\

\noindent\emph{The proof of Theorem \ref{A}.} The proof follows from
Theorem \ref{B} and a simple argument by the Serre duality.\\

Since $L$ is $(n-1)$-ample, by Definition \ref{Def}, for any ample
line bundle $A$,
 there exists a positive number
$m_0=m_0(K_X\ts A^{-1})$ such that when $m>m_0$, we have
$$H^n(X,K_X\ts A^{-1}\ts L^m)=0$$
which is also equivalent to \beq H^0(X,A\ts L^{-m})=0\label{qq}\eeq
by the Serre duality.\\

We argue by contradiction, i.e.  suppose $L$ is $(n-1)$-ample, but
$L$ is not $(n-1)$-positive. In this case,  by Theorem \ref{B}, we
know the dual line bundle $L^{-1}$ must be pseudo-effective. For the
pseudo-effective line bundle $L^{-1}$, it is well-known that, there
exists an ample line bundle $A$ on $X$ such that for every positive
integer $m$ \beq H^0(X,A\ts L^{-m})\neq 0.\label{nonvan}\eeq This
contradicts with (\ref{qq}).\\

 For
readers' convenience, we present a sketched analytical proof of
(\ref{nonvan}) following \cite[Proposition~1.5]{DPS96}. Indeed, for
a fixed very ample line bundle $H$, we choose an ample line bundle
$A$ such that $A\ts K^{-1}_X\ts H^{-n}$ is also ample.
 Hence, there exists a positive rational number $\eps_0$, such that for all $m\geq 0$  and rational number  $\eps\in (0,\eps_0)$,
    $$c_1(L^{-m}\ts A\ts K^{-1}_X\ts H^{-(n+\eps)})$$ lies in the interior of the effective cone of $X$. It implies that
    $L^{-m}\ts A\ts K^{-1}_X\ts H^{-(n+\eps)}$ is linearly equivalent to an effective $\Q$-divisor $D$ plus a
numerically trivial line bundle $T$. Hence \beq \sO_X(L^{-m}\ts
A)\backsimeq \sO_X(K_X \ts H^{(n + \eps)} \ts D \ts T)
\label{iso}\eeq Now, fix a point $x_0\notin D$. Let $\{s_j\} \subset
H^0(X,H)$ be a basis such that all $s_j$ vanish at point $x_0$.  Fix
a local holomorphic basis $e_H$ of $H$ and write $s_j=h_j e_H$. Then
$$h=\frac{1}{\left(\sum_j |h_j|^2\right)^n}$$
is a singular Hermitian metric on $H^n$ with semi-positive curvature
in the sense of current. Moreover, the weight function of $h$ is not
integrable at point $x_0$ and the Lelong number of the curvature
current is $\geq n$. On the other hand, we  can put a singular
metric $h_\eps$ on $H^{\eps}\ts D\ts T$ such that the curvature
current of $h_\eps$ equals $\eps\omega+[D]$ where $\omega$ is a
K\"ahler form in $c_1(H)$ and $[D]$ is the current of integration
over $D$.  Since $x_0\notin D$, the weight function of the singular
metric $hh_\eps$ on $H^{(n + \eps)} \ts D \ts T$ has isolated
singularity at point $x_0$. Moreover,
$$-\sq\p\bp\log(hh_\eps)\geq \eps \omega$$
in the sense of current and the weight function of $hh_\eps$ has
Lelong number $\geq n$ at point $x_0$. By H\"ormander's $L^2$
existence theorem (e.g. \cite[Corollary~3.3]{Dem90}), we know
$$H^0(X, K_X\ts H^{(n + \eps)} \ts D \ts T )\neq 0.$$ By
(\ref{iso}), we obtain (\ref{nonvan}).  The proof of Theorem \ref{A}
is completed.
 \qed

\vskip 2\baselineskip

\noindent \emph{The proof of Proposition \ref{F}}. Let $f:Y\>X$ be the inclusion map. Using the projection formula and the Leray spectral sequence, one has
$$H^i(Y, \mathcal F\ts (f^*L)^{\ts m})=H^i(X, f_*(\mathcal F)\ts L^{\ts m})$$
Hence, if $L\>X$ is $q$-ample, $f^*L\>Y$ is also $q$-ample. On the other hand, since $\dim_\C Y=q+1$
and by Theorem \ref{A}, the $q$-ample line bundle $L|_Y$ is $q$-positive.\qed

\end{document}